\documentclass[letterpaper,11pt]{article}
\usepackage{fullpage,xcolor}
\usepackage{amsmath,amssymb,amsthm,thmtools,mathtools}
\usepackage{cite}
\usepackage{algorithm,algpseudocode}
\usepackage{graphicx,booktabs,multirow}
\usepackage{hyperref}
\usepackage[capitalize,nameinlink]{cleveref}
\usepackage{xjiang}

\allowdisplaybreaks
\mathtoolsset{centercolon}
\setlist[itemize]{leftmargin=1em}
\setlist[enumerate]{leftmargin=1em}
\hypersetup{colorlinks, hypertexnames=false, pageanchor=true,
            linkcolor=blue, citecolor=violet, urlcolor=cyan,
            pdftitle={dcprox},
            }
\makeatletter
\renewcommand{\ALG@name}{Algorithmic~Framework}
\crefname{algorithm}{Framework}{Frameworks}
\makeatother

\newcommand{\primal}{\text{\upshape p}}
\newcommand{\dual}{\text{\upshape d}}
\newcommand{\algo}{{\textsf{\small \upshape DCProx}}}
\newcommand{\mosek}{{\textsf{\small \upshape DC-MOSEK}}}
\newcommand{\fmincon}{{\textsf{\small \upshape fmincon}}}

\theoremstyle{definition}
\newtheorem{assumption}{Assumption}

\theoremstyle{plain}
\newtheorem{theorem}{Theorem}
\newtheorem{prop}{Proposition}

\title{A globally convergent difference-of-convex algorithmic framework and application to log-determinant optimization  problems} 
\author{Chaorui Yao%
\thanks{Department of Electrical and Computer Engineering, UCLA. Email: \textsf{chaorui@ucla.edu}.}
\and Xin Jiang%
\thanks{Department of Industrial and Systems Engineering, Lehigh University. Email: \textsf{xjiang@lehigh.edu}.}
}
\date{June 3, 2023}

\begin{document}
\maketitle

\begin{abstract}
The difference-of-convex algorithm (DCA) is a conceptually simple method for the minimization of (possibly) nonconvex functions that are expressed as the difference of two convex functions.
At each iteration, DCA constructs a global overestimator of the objective and solves the resulting convex subproblem.
Despite its conceptual simplicity, the theoretical understanding and algorithmic framework of DCA needs further investigation.
In this paper, global convergence of DCA at a linear rate is established under an extended Polyak--Łojasiewicz condition.
The proposed condition holds for a class of DC programs with a bounded, closed, and convex constraint set, for which global convergence of DCA cannot be covered by existing analyses.
Moreover, the \algo~computational framework is proposed, in which the DCA subproblems are solved by a primal--dual proximal algorithm with Bregman distances.
With a suitable choice of Bregman distances, \algo~has simple update rules with cheap per-iteration complexity.
As an application, DCA is applied to several fundamental problems in network information theory, for which no existing numerical methods are able to compute the global optimum. 
For these problems, our analysis proves the global convergence of DCA, and more importantly, \algo~solves the DCA subproblems efficiently.
Numerical experiments are conducted to verify the efficiency of \algo.
\end{abstract}

\section{Introduction}

In this paper, we consider the difference-of-convex (DC) programming problems of the form
\begin{equation} \label{eq:prob}
\begin{array}{ll}
  \mini & f(x) = g(x) - h(x) \\
  \st &  x \in \cC,
\end{array}
\end{equation} 
where the functions $g$, $h$ are convex and differentiable, and $\cC$ is a bounded, convex set.
This general problem covers a wide variety of applications in machine learning, information theory, statistics, and other fields \cite{LB16,LP18,SVY18,LNY22}.
A conceptually simple method for solving \eqref{eq:prob} takes the iterations
\begin{equation} \label{eq:dca-iter}
x^{(k+1)} \in \argmin_{x \in \cC} {\big( g(x) - (h(x^{(k)}) + \inprod{\nabla h(x^{(k)})}{x-x^{(k)}}) \big)};
\end{equation} 
\ie, one replaces $h$ with its first-order Taylor expansion at $x^{(k)}$ and then solves the resulting convex subproblem.
Due to the DC structure of \eqref{eq:prob}, the iteration \eqref{eq:dca-iter} will be referred to as the \textit{difference-of-convex algorithm} (DCA).
This paper studies the convergence theory of DCA and numerical methods for the subproblem in~\eqref{eq:dca-iter}, and discusses applications of DCA in a class of matrix optimization problems involving log-determinant functions.

DC programming and DCA have been extensively studied since the late 1990s \cite{PL97,LP05,LB16,LP18}.
Classical convergence analyses for DCA establish asymptotic convergence to a first-order stationary point \cite{PL97,LP18}, and convergence to a global optimum has not been investigated until very recently \cite{WS22a,AdZ23,FFS23}.
However, existing analyses for global convergence only consider the case where $\cC$ is the entire space, and cannot handle closed and bounded constraint sets.
Problems with such more complicated constraint sets arise from network information theory, functional analysis, statistics, \textit{etc.}; see \cite{AJN22,LNY22,WS22a}.
For these problems, only the global minimizer is useful for practical reasons, but unfortunately global convergence of DCA has not been established for them.

Despite the limited global convergence guarantees, numerical methods for the subproblem in~\eqref{eq:dca-iter} are also lacking for the aforementioned applications.
For these problems (and many others), the constraint set $\cC$ is often expensive to project on, which prevents the direct use of first-order methods based on projections and gradient evaluations.
Therefore, the structure of $\cC$ (as well as the objective) should be further exploited, and customized subproblem solvers should be developed for the class of DC programs in the form \eqref{eq:prob}.

In view of the above difficulties, this paper studies the conceptual DC algorithm from both theoretical and algorithmic aspects.
From the theoretical angle, linear global convergence of DCA is established under a new set of conditions.
The proposed assumption extends the recently proposed DC PL condition \cite{FFS23} to account for a bounded, closed, convex constraint set, and is shown to be satisfied in several classes of matrix optimization problems involving log-determinant functions.

In the computational phase, the \algo~algorithmic framework is developed to extend the scope of nonconvex problems numerically solvable by DCA.
In \algo, primal--dual proximal algorithms with Bregman distances are applied to exploit the subproblem structure at each DCA iteration, and solve those subproblems that are not suitable for projected gradient methods.
With a suitable choice of Bregman distances, the inner iterations in \algo~are simple (Bregman) projections and have cheap per-iteration complexity.
As an application, \algo~is applied to several fundamental problems in information theory.
For these problems, \algo~is shown to converge to the global optimum, and has promising numerical performance.

The rest of the paper is organized as follows.
\cref{sec:prelim} reviews the technical background and related work.
The \algo~framework is developed in \cref{sec:algo-dcprox}, and in \cref{sec:algo-conv}, global convergence of DCA is established under a new set of conditions.
\cref{sec:probs} presents several applications, for which DCA converges to the global optimum and customized \algo~has simple update rules with cheap per-iteration complexity.
\cref{sec:results} includes numerical results and \cref{sec:conclusion} concludes the paper.

\section{Background and related work} \label{sec:prelim}

In this section, we discuss some background knowledge on the difference-of-convex algorithm (DCA) and primal--dual first-order proximal algorithms, along with the review on related work.
These two methods are the basic ingredients of the presented algorithmic framework.

We define $[n]=\{1,2,\ldots,n\}$ for $n \in \natint$.
We denote by $\symm^n$ the set of symmetric $n \times n$ matrices, and by $\SV{n}{+}$ the set of symmetric, positive semidefinite (PSD) $n \times n$ matrices.
The inequality $X \succeq 0$ means that $X$ is PSD.
The notation $\inprod{x}{y}=x^Ty$ is used to denote the inner product of vectors $x$ and $y$, and $\|x\|_2 = \inprod{x}{x}^{1/2}$ is the Euclidean norm.
Other norms will be distinguished by subscripts.
The convex conjugate of $f$ is defined as $f^\ast(y) = \sup_x {(\inprod{x}{y} - f(x))}$.

\subsection{The difference-of-convex algorithm (DCA)}

DCA \eqref{eq:dca-iter} can be viewed as a variant of the majorization--minimization (MM) algorithm \cite{SBP17} because at each iteration, it builds a global overestimator of the original problem and solves the resulting simpler subproblem.
The constructed overestimator (known as \textit{surrogate} in MM methods) is able to retain all the information from the convex component and only linearizes the concave portion.
Extensions and variants of \eqref{eq:dca-iter} exist for various kinds of problems.
The current name DCA follows its original paper \cite{PS86}, and the same iteration~\eqref{eq:dca-iter} is also called the \textit{concave--convex procedure} (CCCP) by different researchers; see, \eg, \cite{YR03,SL09,YS22}.
The relation and difference between DCA and CCCP have been discussed in, \eg, \cite{LB16,LP18}, and in this paper, we use the term DCA to refer to the iteration~\eqref{eq:dca-iter} applied to the problem~\eqref{eq:prob}.

Asymptotic convergence of DCA to a first-order stationary point follows naturally from the analysis of MM methods \cite{SBP17}.
Direct analyses of DCA also exist, with different assumptions on the problem; see, \eg, \cite{PL97,SVH05}.
However, global convergence of DCA has not been studied until very recently \cite{WS22a,AdZ23,FFS23}, and is of particular interest here due to the nature of the applications studied in this paper.
Global convergence of DCA can be established when the DC program \eqref{eq:prob} exhibits additional properties.
For example, if the objective in \eqref{eq:prob} is geodesically convex on a manifold $\cM$ and $\cC \subseteq \cM$, DCA \eqref{eq:dca-iter} converges to a global optimum at a linear rate \cite{WS22a}.
In addition, the classical Polyak--Łojasiewicz (PL) condition is also studied and extended for DC programming \cite{AdZ23,FFS23}.
A locally Lipschitz, continuously differentiable function $f$ is said to satisfy the PL condition if there exists $\mu > 0$ such that
\begin{equation} \label{eq:pl}
\mu(f(x) - f(x^\star)) \leq \tfrac{1}{2} \|\nabla f(x)\|_2^2, \quad \text{for all } x \in \dom f,
\end{equation} 
where $x^\star$ is a global minimizer of $f$.
If $\cC = \reals^d$, $g$, $h$ are locally Lipschitz continuous, and $f$ satisfies the PL condition, then DCA \eqref{eq:dca-iter} converges linearly to a global optimum \cite{AdZ23}.
The squared Euclidean distance in the classical PL condition can also be replaced with a Bregman distance (generated by $g^\ast$ or $h^\ast$), yielding the so-called DC PL condition in \cite{FFS23}.
Then linear global convergence is established for DCA if $\cC = \reals^d$, $g$, $h$ are continuously differentiable, and $f$ satisfies the DC PL condition \cite{FFS23}.

\subsection{Primal--dual proximal algorithms with Bregman distances} 

The design of proximal splitting algorithms has recently become an active research area, and proximal methods have been widely applied throughout science and engineering; see \cite{CP11b,CKCH23} for recent surveys.
In general, proximal methods iteratively decompose a large-scale optimization problem into smaller, simpler problems and then solve them separately. 
Thus, proximal methods are suitable for large-scale (or even huge-scale) convex optimization problems that cannot be handled by general-purpose solvers.

An important pillar for proximal algorithms is the (Euclidean) \textit{proximal operator} (or \textit{proximal mapping}) of a closed convex function $f \colon \reals^d \to \reals$, which is defined as
\begin{equation} \label{eq:prox-op}
\prox_f(u) = \argmin_x {\big(f(x) + \tfrac{1}{2} \|x-u\|_2^2 \big)}.
\end{equation}
The minimizer in the definition always exists and is unique for all $u$, as long as $f$ is closed and convex \cite{Moreau65}.
The squared Euclidean distance used in~\eqref{eq:prox-op} can also be replaced with a generalized distance~$d$, of which an important example is the Bregman distance. 
Let $\phi$ be a convex function with a domain that has a nonempty interior, and assume $\phi$ is continuous on $\dom \phi$ and continuously differentiable on $\intr(\dom \phi)$. The \textit{Bregman distance} \cite{CZ97} generated by the \textit{kernel function} $\phi$ is 
\[
d(x,y) = \phi(x) - \phi(y) - \inprod{\nabla \phi(y)}{x-y},  
\]
with $\dom d = \dom \phi \times \intr (\dom \phi)$. 
The corresponding \textit{Bregman proximal operator} is defined as:
\[
\prox_f^\phi (u,a) = \argmin_x {\big(f(x) + \inprod{a}{x} + d(x,u) \big)}.
\]
It is assumed that the minimizer in the definition is unique and lies in $\intr(\dom \phi)$ for all $a$ and all $u \in \intr (\dom \phi)$.
In all the applications studied in this paper, the Bregman proximal operator is used with specific combinations of $f$ and $\phi$, and its existence and uniqueness follow directly from a closed-form solution.
The use of Bregman distances can make the corresponding proximal operator easier to compute, and/or help build a more accurate local optimization model around the current iterate of a proximal algorithm \cite{JV22b}.
We refer interested readers to \cite{BBT17,BSTV18,CLTY23,CV18} for more applications of Bregman proximal methods.

In this paper, proximal algorithms are applied to convex problems of the form
\begin{equation} \label{eq:prob-pdhg}
\mini \quad F(u) + G(\cA u),
\end{equation} 
with optimization variable $u$, where $F$, $G$ are closed convex functions and $\cA$ is a linear mapping. A well-known proximal method for solving this type of problem is the Bregman primal--dual hybrid-gradient (PDHG) algorithm \cite{CP11a,CP16b}:
\begin{subequations} \label{eq:breg-pdhg}
\begin{align}
u^{(t+1)} &= \prox_{\tau F}^{\phi_\primal} (u^{(t)}, \tau \cA^\ast v^{(t)}) \label{eq:breg-pdhg-x} \\ 
v^{(t+1)} &= \prox_{\sigma G^\ast}^{\phi_\dual} (v^{(t)}, -\sigma \cA(2u^{(t+1)}-u^{(t)})), \label{eq:breg-pdhg-z}
\end{align}
\end{subequations} 
where $\cA^\ast$ is the adjoint operator of $\cA$, and $\phi_\primal$ and $\phi_\dual$ are two kernel functions in the primal and dual spaces, respectively.
It is assumed (without loss of generality) that $\phi_\primal$ and $\phi_\dual$ are 1-strongly convex with respect to norms $\|\cdot\|_\primal$ and $\|\cdot\|_\dual$:
\[
d_\primal (u,u^\prime) \geq \tfrac{1}{2} \|u-u^\prime\|_\primal^2, \qquad d_\dual (v,v^\prime) \geq \tfrac{1}{2} \|v-v^\prime\|_\dual^2
\]
for all $(u,u^\prime) \in \dom d_\primal$ and $(v,v^\prime) \in \dom d_\dual$.
The stepsizes $\sigma$, $\tau$ must satisfy $\sigma \tau \|\cA\|^2 \leq 1$. Here the operator norm $\|\cA\|$ is defined as 
\[
\|\cA\| = \sup_{u \neq 0, v \neq 0} \frac{\inprod{v}{Au}}{\|v\|_\dual \|u\|_\primal} = \sup_{u \neq 0} \frac{\|\cA u\|_{\dual,\ast}}{\|u\|_\primal} = \sup_{v \neq 0} \frac{\|\cA^\ast v\|_{\primal,\ast}}{\|v\|_\dual},
\]
where $\|\cdot\|_{\primal,\ast}$ and $\|\cdot\|_{\dual,\ast}$ are the dual norm of $\|\cdot\|_\primal$ and $\|\cdot\|_\dual$, respectively.
Line search techniques also exist for choosing adaptive stepsizes \cite{MP18,JV22a,JV22b}.
With suitable stepsizes, the iterates $(u^{(k)}, v^{(k)})$ converge to a primal and dual optimum of the problem~\eqref{eq:prob-pdhg}, and the averaged iterates converge at a rate of $O(1/k)$ \cite{Condat13,CP16b,JV22b}.
Successful applications of Bregman PDHG range from signal processing \cite{CV18}, optimal transport \cite{CC22} to sparse semidefinite programming \cite{JV22a}.

\section{\textsf{DCProx} algorithmic framework and global convergence of DCA} \label{sec:algo}

We now present our theoretical and algorithmic contributions to DCA.
When DCA is applied to the problem \eqref{eq:prob} (with a nontrivial constraint set~$\cC$), global convergence cannot be obtained directly from existing results, and the subproblem \eqref{eq:dca-iter} cannot be solved efficiently by first-order methods based on projections and gradient evaluations.
However, when $\cC$ exhibits special structure (\eg, it is written as the intersection of simple convex sets), primal--dual proximal methods can be used to split the nontrivial set constraints and each (primal/dual) update becomes a simple (Bregman) projection.
In addition, the recently proposed DC PL condition~\cite{FFS23} is extended to account for a closed and bounded set~$\cC$, and global convergence is established for DCA when applied to~\eqref{eq:prob}.

\subsection{The \textsf{DCProx} framework} \label{sec:algo-dcprox}

Very often in applications, the constraint set $\cC$ is nontrivial, and projection onto $\cC$ is expensive.
Thus, most projection-and-gradient-based first-order methods are not efficient in solving the subproblem in~\eqref{eq:dca-iter}, and special structure of $\cC$ should be leveraged to design a customized subproblem solver.
In many applications, the set $\cC$ can be naturally written as the intersection of two (or several) simple convex sets, and projection onto each component set is much easier than projection on $\cC$.
(This is true for all the applications in \cref{sec:probs}, and in particular, \cref{sec:prob-3} provides an example in which~$\cC$ consists of an arbitrary number of component sets.)
This motivates the use of primal--dual proximal algorithms (\eg, PDHG \cite{CP11a,CP16b}) as the subproblem solver.
For example, if $\cC = \cC_1 \cap \cC_2$, where $\cC_1$ and $\cC_2$ are convex sets, the subproblem~\eqref{eq:dca-iter} can be reformulated in the form of \eqref{eq:prob-pdhg} with
\[
F(u) = g(u) - \inprod{\nabla h(x^{(k)})}{u} + \delta_{\cC_1} (u), \qquad G(v) = \delta_{\cC_2} (v), \qquad \cA = \mathrm{Id},
\]
where $\delta_{\cC_i}$ is the indicator function of $\cC_i$ and $\mathrm{Id}$ is the identity operator.
Therefore, with Bregman PDHG~\eqref{eq:breg-pdhg} as the subproblem solver, the proposed \algo~algorithmic framework is summarized in \cref{algo}.
(The term \textsf{\small Prox} indicates that proximal methods exploit a splitting of $\cC$.)
\begin{algorithm}[t]
\caption{\algo}
\label{algo}
\begin{algorithmic}[1]
\For{$k=1,2,\ldots,$}
  \For{$t=1,2,\ldots,$}
    \State Compute the Bregman PDHG updates:
      \begin{subequations} \label{eq:algo}
        \begin{align}
          u^{(t+1)} &\leftarrow \argmin_{u \in \cC_1} {g(u) + \inprod{\cA^\ast v^{(t)} - \nabla h(x^{(k)})}{u} + \tfrac{1}{\tau} d_\primal (u,u^{(t)}) } \label{eq:algo-u} \\
          v^{(t+1)} &\leftarrow \argmin_{v} {\delta_{\cC_2}^\ast (v) - \inprod{\cA (2u^{(t+1)} - u^{(t)})}{v} + \tfrac{1}{\sigma} d_\dual (v,v^{(t)}) }. \label{eq:algo-v}
        \end{align}
      \end{subequations} 
  \EndFor

  \State Update DCA iterate $x^{(k+1)} \leftarrow \uhat$, where $(\uhat,\vhat)$ is the limit point of Bregman PDHG. 
\EndFor
\end{algorithmic}
\end{algorithm}
It is called an \textit{algorithmic framework} rather than an algorithm because in principle, $\phi_\primal$ and $\phi_\dual$ could be any valid kernel functions.
The choice of Bregman kernels largely depends on problem structure, and a suitable choice of kernels can significantly improve the efficiency of \algo.
The distance $d_\dual$ in the dual space is often chosen as the squared Euclidean distance ($\phi_\dual = \tfrac{1}{2} \|\cdot\|_2^2$), and then the dual update~\eqref{eq:algo-v} simply involves the Euclidean projection on~$\cC_2$.
Moreover, when we choose $\phi_\primal = g$, the primal update \eqref{eq:algo-u} reduces to
\begin{align*}
u^{(t+1)} &= \argmin_{u \in \cC_1} {(1+\tfrac{1}{\tau}) g(u) + \big\langle {\tfrac{\tau}{1+\tau} (\cA^\ast v^{(t)} - \nabla h(x^{(k)})) - \tfrac{1}{\tau} \nabla g(u^{(t)})}, {u} \big\rangle } \\
&= \argmin_{u \in \cC_1} d_g (u, \utilde), \quad \text{where } \nabla g(\utilde) = \tfrac{\tau}{1+\tau} (\cA^\ast v^{(t)} - \nabla h(x^{(k)})) - \tfrac{1}{\tau} \nabla g(u^{(t)}).
\end{align*}
With this choice of the kernel, the primal update is simply a Bregman projection on $\cC_1$.
(Note that~$x^{(k)}$ is the $k$th DCA iterate while $(u^{(t)}, v^{(t)})$ is the $t$th PDHG iterate.
The vector $\utilde$ is used only for demonstration purpose, and is not computed explicitly in practice.)

\subsection{Global convergence for DCA} \label{sec:algo-conv}

While the presented \algo~framework makes DCA tractable for a wider range of DC problems (with a nontrivial convex constraint set), we show that for these problems, DCA converges to a global optimum at a linear rate.
(Note that the presented analysis focuses on the DCA iterations~\eqref{eq:dca-iter}, and is independent of the subproblem solver.)
All the assumptions for \eqref{eq:prob} are listed below, and global convergence of DCA for this kind of problems are not addressed by existing analyses \cite{AdZ23,FFS23,WS22a}.
\begin{assumption} \label{assume}
The following assumptions are made for the problem \eqref{eq:prob}.
\begin{enumerate}
\item The functions $g$, $h$ are closed, convex, and continuously differentiable on their respective domains.

\item The constraint set $\cC$ is bounded, closed, and convex, with $\cC \subseteq \dom g \cap \dom h$.

\item The problem \eqref{eq:prob} has a global minimizer $x^\star$ with a finite optimal value $f^\star = f(x^\star)$.

\item The DCA iteration~\eqref{eq:dca-iter} is solvable for any point $x^{(k)} \in \cC$.

\item There exists $\mu > 0$ such that for some $r>0$,
\begin{equation} \label{eq:dc-pl}
\mu (f(x) - f^\star) \leq d_{h^\ast} (\nabla g(x) + y, \nabla h(x)), \quad \text{for all } x \in \cC, \ y \in N_\cC(x) \ \text{and} \ \|y\|_2 \leq r
\end{equation} 
where $N_\cC(x)$ is the normal cone of $\cC$ at $x$.
\end{enumerate}
\end{assumption}
The differentiability of $g$ and $h$ implies that $\dom f$ is an open set, and existing analyses \cite{AdZ23,FFS23} that (implicitly) assume $\cC=\reals^d$ cannot handle a closed and bounded constraint set.
\cref{assume}.4 is an extension of the recently proposed DC PL condition \cite[Definition~1]{FFS23} and the classical PL condition.
(A similar condition is proposed in \cite{BBC+19} for the mirror descent method.)
In particular, \eqref{eq:dc-pl} replaces the squared Euclidean distance in the PL condition~\eqref{eq:pl} with a Bregman distance, and carefully handles the closedness and boundedness of $\cC$ that is often ignored in the DCA literature.
The constant $r>0$ is needed to bound $\|y\|_2$ when $y$ is on the boundary of $\cC$, and is used only for theoretical analysis.
\cref{assume}.4 might (implicitly) impose some requirements on the constant~$r$, but its actual value does not need to be known in practice.
With \cref{assume}, we can prove a linear convergence rate on the function value.
\begin{theorem} \label{thm:dca-conv}
With \cref{assume}, the DCA iteration \eqref{eq:dca-iter} satisfies
\[
f(x^{(k+1)}) - f^\star \leq \frac{1}{1+\mu} (f(x^{(k)}) - f^\star).
\] 
\end{theorem}
\begin{proof}
From~\cref{assume}, $\|\nabla g(x)\|_2$ and $\|\nabla h(x)\|_2$ are bounded for all $x \in \cC$, and then it is assumed that the constant $r$ satisfies $r > \|\nabla g(x)\|_2 + \|\nabla h(x)\|_2$ for all $x \in \cC$.
The optimality condition of DCA implies that there exists $\yhat \in N_\cC(x^{(k+1)})$ and $\|\yhat\|_2 \leq r$ such that $\nabla h(x^{(k)}) = \nabla g(x^{(k+1)}) + \yhat$. Then it follows from \eqref{eq:dc-pl} that
\begin{align*}
\MoveEqLeft[0.2] \mu (f(x^{(k+1)}) - f^\star) \\ 
&\leq d_{h^\ast} (\nabla g(x^{(k+1)}) + \yhat, \nabla h(x^{(k+1)})) \\
&= d_{h^\ast} (\nabla h(x^{(k)}), \nabla h(x^{(k+1)})) \\ 
&= h^\ast (\nabla h(x^{(k)})) - h^\ast (\nabla h(x^{(k+1)})) - \inprod{x^{(k+1)}}{\nabla h(x^{(k)}) - \nabla h(x^{(k+1)})} \\ 
&= \inprod{\nabla h(x^{(k)})}{x^{(k)}} - h(x^{(k)}) - \inprod{\nabla h(x^{(k+1)})}{x^{(k+1)}} + h(x^{(k+1)}) \\ 
&\phantom{=} \mathrel{-} \inprod{x^{(k+1)}}{\nabla h(x^{(k)})-\nabla h(x^{(k)})} \\ 
&= d_h (x^{(k+1)}, x^{(k)}) \\ 
&= f(x^{(k)}) - f(x^{(k+1)}) - d_g (x^{(k)}, x^{(k+1)}) \\ 
&\leq (f(x^{(k)}) - f^\star) - (f(x^{(k+1)}) - f^\star).
\end{align*}
Thus, the desired inequality follows.
\end{proof}

\section{Application to log-determinant optimization problems} \label{sec:probs}

In this section, we apply \algo~to several classes of problems arising from network information theory \cite{WSS06,GN14,AJN22,LNY22}.
Each of these problems is shown to have a unique global solution \cite{LNY22}, but so far we are not aware of any optimization algorithm that is guaranteed to find the optimum.
In comparison, with \cref{thm:dca-conv} we show that DCA converges to the unique global optimum for all the problems of interest.
Moreover, for all the studied applications, projected gradient methods might not be applicable to the subproblems in \eqref{eq:dca-iter} due to a nontrivial $\cC$.
In comparison, \algo~is able to solve these subproblems by splitting the nontrivial constraint set $\cC$, and in particular, the use of a suitable Bregman distance further improves the efficiency of \algo.

\subsection{Two receiver vector Gaussian broadcast channel with private messages} \label{sec:prob-1}

Consider the optimization problem
\begin{equation} \label{eq:1-prob}
\begin{array}{ll}
  \mini & -\log\det (X+\Sigma_1) + \lambda \log\det (X+\Sigma_2) \\
  \st & 0 \preceq X \preceq C,
\end{array}
\end{equation} 
where the variable is $X \in \symm^n$, and the coefficients are $\Sigma_1,\Sigma_2 \in \SV{n}{++}$, $C \in \SV{n}{+}$, and $\lambda > 1$.
This problem arises from various applications in network information theory, wireless communication, \textit{etc.}; see, \eg, \cite{WSS06,LV07}.
In certain applications, the unique global optimum of the problem evaluates the capacity region of two receiver vector Gaussian broadcast channels with private messages \cite{LNY22}.
Due to the nature of its origin, only the global solution of~\eqref{eq:1-prob} is useful for the underlying application, and other first-order stationary points are of limited interest.

\subsubsection{DCA and its global convergence}

The problem \eqref{eq:1-prob} is a DC program in the form of \eqref{eq:prob} with
\[
g(X) = -\log \det (X+\Sigma_1), \qquad h(X) = -\lambda \log \det (X+\Sigma_2), \qquad \cC = \{X \in \symm^n \mid 0 \preceq X \preceq C\}.
\]
In this problem, the functions $g$, $h$ are strongly convex on the bounded set $\cC$, and have an Lipschitz continuous gradient on $\intr \cC$.
(But they are not strongly convex nor $L$-smooth on their natural domains.)
This fact helps to establish the inequality \eqref{eq:dc-pl} for \eqref{eq:1-prob}, and then global convergence of DCA is provided in \cref{prop:1-conv}.
\begin{prop} \label{prop:1-conv}
Problem \eqref{eq:1-prob} satisfies \cref{assume}. 
DCA for \eqref{eq:1-prob} converges to its minimizer at a linear rate.
\end{prop}
\begin{proof}
It suffices to show that the problem \eqref{eq:1-prob} satisfies the inequality \eqref{eq:dc-pl}.
Within the bounded, convex set $\cC = \{X \in \symm^n \mid 0 \preceq X \preceq C\}$, the function $g(X) = -\log \det (X+\Sigma_1)$ has a Lipschitz continuous gradient and $h(X) = -\lambda \log \det (X+\Sigma_2)$ is strongly convex. Thus, there exists two positive constants $\alpha, \beta > 0$ such that for all $X \in \cC$,
\begin{align*}
  g(X) - g(X^\star) &\leq \tfrac{\alpha}{2} \|X-X^\star\|_F^2 \\
  h(X) - h(X^\star) &\geq \tfrac{\beta}{2} \|X-X^\star\|_F^2,
\end{align*}
where $X^\star$ is the global optimum for \eqref{eq:1-prob}.
Thus it follows that
\begin{equation} \label{eq:1-conv-prf-1}
  f(X) - f(X^\star) \leq \tfrac{\alpha-\beta}{2} \|X-X^\star\|_F^2
\end{equation} 
for all $X \in \cC$.
(It implies from the optimality of $X^\star$ that $\alpha > \beta$.)

On the other hand, the convex conjugate $h^\ast$ is also strongly convex on the bounded, convex set~$\cC$. Thus, there exists some $\gamma > 0$ such that $d_{h^\ast}(U,V) \geq \tfrac{\gamma}{2} \|U-V\|_F^2$; see, \eg, \cite{JV22b}. Hence, the right-hand side of \eqref{eq:dc-pl} is bounded below by
\begin{align}
d_{h^\ast} (\nabla g(X) + Y, \nabla h(X)) &\geq \tfrac{\gamma}{2} \|\nabla g(X) + Y - \nabla h(X)\|_F^2 \nonumber \\
&= \tfrac{\gamma}{2} \|\nabla g(X) + Y - \nabla h(X) - (\nabla g(X^\star) + Y^\star - \nabla h(X^\star))\|_F^2 \nonumber \\
&\geq \tfrac{\gamma}{2} \big| \|\nabla f(X) - \nabla f(X^\star) + (Y-Y^\star)\|_F^2, \label{eq:1-conv-prf-2}
\end{align}
where $Y^\star \in N_\cC(X^\star)$ satisfies $\nabla g(X^\star) + Y^\star = \nabla h(X^\star)$.
Following the same argument in the first part, we can see that there exists $\kappa_\mathrm{max} > \kappa_\mathrm{min} > 0$ such that $\kappa_\mathrm{min} \|X-X^\star\|_F \leq \|\nabla f(X) - \nabla f(X^\star)\|_F \leq \kappa_\mathrm{max} \|X-X^\star\|_F$ for all $X \in \cC$.
In addition, for all $Y \in N_\cC(X)$ and $\|Y\|_F \leq r$, there exists $\xi_\mathrm{max} > \xi_\mathrm{min} > 0$ such that $\xi_\mathrm{min} \leq \|Y-Y^\star\|_F \leq \xi_\mathrm{max}$.
Thus, it follows that
\begin{align}
\MoveEqLeft[0.2] \|\nabla f(X) - \nabla f(X^\star) + (Y-Y^\star)\|_F \nonumber \\
&\geq \big|\|\nabla f(X)-\nabla f(X^\star)\|_F - \|Y-Y^\star\|_F \big| \nonumber \\ 
&\geq \max\{\kappa_\mathrm{min} \|X-X^\star\|_F - \xi_\mathrm{max}, \xi_\mathrm{min} - \kappa_\mathrm{max} \|X-X^\star\|_F\}.
\end{align}
Clearly, there exists $\eta > 0$ such that
\begin{equation} \label{eq:1-conv-prf-4}
  \eta \|X-X^\star\|_F \leq \xi_\mathrm{min} - \kappa_\mathrm{max} \|X-X^\star\|_F.
\end{equation}
Combining \eqref{eq:1-conv-prf-1}--\eqref{eq:1-conv-prf-4}, we conclude that there exists $\mu = (\alpha-\beta)/ (\eta^2 \gamma) > 0$ such that
\[
  \mu (f(X) - f(X^\star)) \leq d_{h^\ast} (\nabla g(X) + Y, \nabla h(X)). \qedhere
\]
\end{proof}

\subsubsection{\algo: Bregman PDHG as the subproblem solver}

The $k$th DCA iteration \eqref{eq:dca-iter} for problem \eqref{eq:1-prob} involves the convex subproblem
\begin{equation} \label{eq:1-subprob}
\begin{array}{ll}
\mini & -\log \det (U+\Sigma_1) + \lambda \inprod{(X^{(k)}+\Sigma_2)^{-1}}{U} \\
\st & 0 \preceq U \preceq C
\end{array}
\end{equation}
with optimization variable $U \in \symm^n$.
Projection onto the constraint set~$\cC$ is much more expensive than a PSD projection, which prevents the direct use of first-order methods based on projections and gradient evaluations.
Nevertheless, the constraint set can be written as $\cC = \SV{n}{+} \cap \cC_0$ with $\cC_0 = \{U \in \symm^n \mid U \preceq C\}$, and projection on each component set ($\SV{n}{+}$ and $\cC_0$) is much easier.
This observation motivates the use of primal--dual proximal algorithms, and when \cref{algo} is applied, the inner iterations~\eqref{eq:algo} (with $\phi_\primal=g$ and $\phi_\dual=\tfrac{1}{2}\|\cdot\|_2^2$) reduce to
\begin{align*}
U^{(t+1)} &= \argmin_{U \in \SV{n}{+}} {-(1+\tfrac{1}{\tau}) \log \det (U+\Sigma_1) + \inprod{V^{(t)}+\lambda (X^{(k)}+\Sigma_2)^{-1}+\frac{1}{\tau} (U^{(t)}+\Sigma_1)^{-1}}{U} } \\
V^{(t+1)} &= \prox_{\sigma \delta_{\cC_0}^\ast}
  (V^{(t)} + \sigma(2U^{(t+1)} - U^{(t)})),
\end{align*}

\paragraph{Dual update}
Recall $\prox_{\delta_{\cC_0}} = \proj_{\cC_0}$ is the projection operator, and the $V$-update reduces to
\[
V^{(t+1)} = \proj_{\SV{n}{+}} \big(\sigma C-V^{(t)}-\sigma (2U^{(t+1)} - U^{(t)}) \big) - \big(\sigma C - V^{(t)}-\sigma (2U^{(t+1)} - U^{(t)}) \big);
\]
\ie, given the eigen-decomposition of $\Vhat = \sigma C - V^{(t)} - \sigma (2U^{(t+1)} - U^{(t)})$, the next iterate $V^{(t+1)}$ simply sets the positive eigenvalues of $\Vhat$ to zero and takes the absolute values of the negative ones. 

\paragraph{Primal update}
If the squared Euclidean distance is used in the $U$-update, the computation of the proximal mapping needs an eigen-decomposition and solving a number of $n$ quadratic equations (see, \eg, \cite[\S6.5]{BPC+11}).
In comparison, with $\phi_\primal(X)=g(X)$, we can avoid solving a large number of quadratic equations, and the $U$-update has a closed-form expression given by \cref{prop:1-logdet-sol}.
\begin{prop} \label{prop:1-logdet-sol}
Consider the optimization problem
\[
\begin{array}{ll}
\mini & -\log\det (X+\Sigma) + \inprod{A}{X} \\
\st & X \succeq 0
\end{array}
\]
with variable $X \in \symm^n$ and data matrices $A, \Sigma \in \SV{n}{++}$. 
Suppose $\Sigma^{1/2} A \Sigma^{1/2}=Q \Lambda Q^T$ is the eigen-decomposition, with eigenvalues $\lambda_1 \geq \lambda_2 \geq \cdots \geq \lambda_n > 0$. 
Then the optimal solution is
\[
X^\star = \Sigma^{1/2} Q \psi(\Lambda) Q^T \Sigma^{1/2}.
\]
Here the scalar function $\psi$ is defined as $\psi(\gamma) = \max\{(1-\gamma)/\gamma, 0\}$, and $\psi(\Lambda)$ is a diagonal matrix with diagonal elements $\psi(\Lambda)_{ii} = \psi(\lambda_i)$.
\end{prop}
\begin{proof}
With a change of variables $Y=Q^T \Sigma^{-1/2} X \Sigma^{-1/2} Q$, the objective function can be equivalently written as
\begin{align*}
f(X) &= -\log \det \big(\Sigma^{1/2} (I+\Sigma^{-1/2} X \Sigma^{-1/2}) \Sigma^{1/2} \big)
  + \inprod{\Sigma^{1/2} A \Sigma^{1/2}}{\Sigma^{-1/2} X \Sigma^{-1/2}} \\
&= -\log \det \Sigma - \log \det (I+\Sigma^{-1/2} X \Sigma^{-1/2}) + \inprod{\Sigma^{1/2} A \Sigma^{1/2}}{\Sigma^{-1/2} X \Sigma^{-1/2}} \\
&= -\log \det \Sigma - \log \det (I+Y) + \inprod{\Lambda}{Y}.
\end{align*}
We argue that the optimal $Y$ is a diagonal matrix. To see this, suppose $Y$ is diagonal. Then it follows that $\inprod{\Lambda}{Y} \leq \inprod{\Lambda}{VYV^T}$, for any diagonal matrix $\Lambda$ with nonnegative diagonal elements and any orthogonal matrix $V$.
Hence the optimal $Y$ has to diagonal, and $f$ further reduces to
\[
f(X) = -\log \det \Sigma + \sum_{i=1}^n (\lambda_i Y_{ii} - \log (1+Y_{ii})).
\] 

Recall that $A, \Sigma \in \SV{n}{++}$, and it follows that $\gamma_i > 0$ for all $i=1,\ldots,n$.
Define the scalar function $h(\gamma) = \lambda \gamma - \log(1+\gamma)$ with domain $\reals_+$.
If $\lambda > 1$, the function $h$ is monotone increasing on its domain, and the minimizer of $h$ is taken at $0$.
If $0 < \lambda \leq 1$, the function $h$ takes its minimum at $\gamma=(1-\lambda)/\lambda$.
Hence, the optimal $Y^\star$ is a diagonal matrix with diagonal entries
\[
Y_{ii}^\star = \psi(\lambda_i) = \max \Big\{\frac{1-\lambda_i}{\lambda_i}, 0\Big\},
\]
and then the optimal solution $X^\star$ is given by
\[
X^\star = \Sigma^{1/2} Q \begin{bmatrix}
  \psi(\lambda_1) & & & \\ & \psi(\lambda_2) & & \\ & & \ddots & \\ & & & \psi(\lambda_n)
\end{bmatrix} Q^T \Sigma^{1/2}. \qedhere
\]
\end{proof}

\subsection{Gaussian broadcast channel with common messages} \label{sec:prob-2}

The second nonconvex DC program of interest is a generalization of \eqref{eq:1-prob}, and reads as
\begin{equation} \label{eq:2-prob}
  \begin{array}{ll}
    \mini & - \beta \log\det (X+Y+\Sigma_2) 
      + \alpha \log\det (X+Y+\Sigma_1) \\
    & - \log\det (X+\Sigma_1) + \lambda \log\det(X+\Sigma_2) \\
    \st & X + Y \preceq C, \; X \succeq 0, \; Y \succeq 0,
  \end{array}
\end{equation} 
where the optimization variables are $X, Y \in \symm^n$, and the coefficients are $\Sigma_1, \Sigma_2 \in \SV{n}{++}$, $C \in \SV{n}{+}$, $\alpha \in [0,1]$, $\beta > 0$, and $\lambda > 1$. 
(This choice of parameters is included in \cite{GN14}, and other choices in \cite{GN14} give an optimization problem that is easier to solve than the presented one.)
The unique solution of \eqref{eq:2-prob} evaluates the capacity region of two receiver vector Gaussian broadcast channels with common messages and private messages; see \cite{GN14} and references therein.

The problem~\eqref{eq:2-prob} can also be written in the form of \eqref{eq:prob}, with convex functions $g, h \colon \symm^n \times \symm^n \to \reals$
\begin{align*}
    g(X,Y) &= - \beta \log\det (X+Y+\Sigma_2) - \log\det (X+\Sigma_1), \\
    h(X,Y) &= -\alpha \log\det (X+Y+\Sigma_1) - \lambda \log\det(X+\Sigma_2),
\end{align*}
and the constraint set $\cC = \{(X,Y) \mid X \succeq 0, Y \succeq 0, X+Y \preceq C\}$.
Following the same ideas as in \cref{prop:1-conv}, we show that DCA applied to \eqref{eq:2-prob} converges to its unique global optimum.
\begin{prop} \label{prop:2-conv}
Problem \eqref{eq:2-prob} satisfies \cref{assume}. DCA for \eqref{eq:2-prob} converges to its minimizer at a linear rate.
\end{prop} 

In the remainder of this subsection, we focus on the convex subproblem at $k$th DCA iteration:
\begin{equation} \label{eq:2-subprob}
\begin{array}{ll}
\mini & -\beta \log\det (U + \Utilde + \Sigma_2)
  +\alpha \inprod{(X^{(k)} + Y^{(k)} + \Sigma_1)^{-1}}{U + \Utilde} \\
& -\log\det (U + \Sigma_1) + \lambda \inprod{(X^{(k)} + \Sigma_2)^{-1}}{U} \\
\st & U + \Utilde \preceq C, \; U \succeq 0, \; \Utilde \succeq 0
\end{array}
\end{equation} 
with variables $U, \Utilde \in \symm^n$.
We first make a change of variables $W := U + \Utilde$, and then reformulate the problem~\eqref{eq:2-subprob} in the form of~\eqref{eq:prob-pdhg} with
\begin{align*}
F(U,W) &= -\beta \log\det (W + \Sigma_2)
  + \alpha \inprod{(X^{(k)} + Y^{(k)} + \Sigma_1)^{-1}}{W} \\
&\phantom{=} \mathrel{-} \log\det (U + \Sigma_1)
  + \lambda \inprod{(X^{(k)} + \Sigma_2)^{-1}}{U}
  + \delta_{\SV{n}{+}}(U) + \delta_{\cC_0} (W) \\
G &= \delta_{\SV{n}{+}}, \quad 
\cC_0 = \{W \in \symm^n \mid W \preceq C\},
\end{align*}
and the linear mapping $\cA \colon \symm^n \times \symm^n \to \symm^n$ is defined as $\cA(U,W) = W-U$.
When \cref{algo} is applied to this splitting strategy, the inner iterations \eqref{eq:algo} reduce to
\begin{align*}
U^{(t+1)} &= \argmin_{U \in \SV{n}{+}} \big(-\log\det (U+\Sigma_1)
  + \inprod{-V^{(t)} + \lambda (X^{(k)}+\Sigma_2)^{-1}}{U} 
  + \tfrac{1}{\tau} d_1 (U,U^{(t)}) \big) \\
W^{(t+1)} &= \argmin_{W \in \cC} \big(-\beta \log\det (W+\Sigma_2)
  +\inprod{V^{(t)} + \alpha (X^{(k)}+Y^{(k)}+\Sigma_1)^{-1}}{W}
  + \tfrac{1}{\tau} d_2 (W,W^{(t)}) \big) \\
V^{(t+1)} &= \prox_{\sigma \delta^\ast_{\SV{n}{+}}}
  \big(V^{(t)} + \sigma (2(W^{(t+1)}-U^{(t+1)})-(W^{(t)}-U^{(t)})) \big).
\end{align*}
The $V$-update is simple and involves an eigen-decomposition.
In the $U$-update, we use the Bregman distance generated by $\phi_1(V) = -\log\det(V+\Sigma_1)$,
and the closed-form update rule is given by \cref{prop:1-logdet-sol}.
A similar result exists for the $W$-update if we choose the kernel function $\phi_2(W)=-\log\det(W+\Sigma_2)$ with $\dom \phi_2 = \{W \in \symm^n \mid W + \Sigma_2 \succ 0\}$.
Then $W$-update minimizes 
\[
\ftilde (W) = -(\beta +\tfrac{1}{\tau}) \log \det (W+\Sigma_2)
+ \inprod{V^{(t)}+\alpha (X^{(k)}+Y^{(k)}+\Sigma_1)^{-1}
+ \tfrac{1}{\tau} (W^{(t)}+\Sigma_2)^{-1}}{W} 
\]
subject to $W \preceq C$, and has a closed-form expression given in \cref{prop:2-logdet-sol}.

\begin{prop} \label{prop:2-logdet-sol}
Consider the optimization problem
\[
\begin{array}{ll}
\mini & -\log\det (X+\Sigma) + \inprod{A}{X} \\
\st &  X \preceq C
\end{array}
\]
with variable $X \in \symm^n$ and data $A, \Sigma \in \SV{n}{++}$, $C \in \SV{n}{+}$.
Suppose $(C+\Sigma)^{1/2} A (C+\Sigma)^{1/2}=Q \Lambda Q^T$ is the eigen-decomposition, with eigenvalues $\lambda_1 \geq \lambda_2 \geq \cdots \geq \lambda_n > 0$. 
Then the optimal solution is
\[
X^\star = (C+\Sigma)^{1/2} Q \psi(\Lambda) Q^T (C+\Sigma)^{1/2} - \Sigma.
\]
Here the scalar function $\psi$ is defined as $\psi(\gamma) = \min\{1/\gamma, 1\}$, and $\psi(\Lambda)$ is a diagonal matrix with diagonal elements $\psi(\Lambda)_{ii} = \psi(\lambda_i)$.
\end{prop}
\begin{proof}
With a change of variables $Z=X+\Sigma, Y=Q^T (C+\Sigma)^{-1/2}Z(K+\Sigma)^{-1/2} Q$, the objective function can be equivalently written as
\[
f(X) = -\log \det (C+\Sigma) - \inprod{A}{\Sigma} - \log \det Y + \inprod{\Lambda}{Y},
\] 
and the constraint is equivalent to $Y \preceq I$. Following the argument in the proof of \cref{prop:1-logdet-sol}, we can show that the optimal $Y$ is a diagonal matrix. Then $f$ further reduces to
\[
f(X) = -\log \det (C+\Sigma) - \inprod{A}{\Sigma} + \sum_{i=1}^n \lambda_i Y_{ii} - Y_{ii}.
\] 

Define the scalar function $h(\gamma) = \lambda \gamma - \log \gamma$ and parameter $\lambda \in \reals_{++}$.
The function $h$ is decreasing on $(0,1/\lambda)$ and increasing on $(1/\lambda,+\infty)$.
Hence, the optimal $Y^\star$ is a diagonal matrix with diagonal entries
\[
Y_{ii}^\star = \psi(\lambda_i) = \max \Big\{\frac{1-\lambda_i}{\lambda_i}, 0\Big\},
\] 
and the optimal solution $X^\star$ is given by
\[
X^\star = (C+\Sigma)^{1/2} Q \begin{bmatrix}
  \psi(\lambda_1) & & & \\ & \psi(\lambda_2) & & \\ & & \ddots & \\ & & & \psi(\lambda_n)
\end{bmatrix} Q^T (C+\Sigma)^{1/2} - \Sigma. \qedhere
\]
\end{proof}

\subsection{Generalized Brascamp--Lieb inequality} \label{sec:prob-3}

The third problem arises from theoretical computer science and functional analysis \cite{BL76,BCCT08,AJN22}
\begin{equation} \label{eq:3-prob}
\begin{array}{ll}
\mini & -\sum\limits_{i=1}^p \beta_i \log \det X_i + \sum\limits_{j=1}^q \alpha_j \log \det \big(\sum\limits_{i=1}^p A_{ij} X_i A_{ij}^T + \rho I_{m_j} \big) \\ 
\st & 0 \preceq X_i \preceq C_i, \quad i \in [p],
\end{array}
\end{equation} 
where the optimization variables are $X_i \in \symm^{n_i}$, $i \in [p]$, and the coefficients are $C_i \in \SV{n_i}{+}$ for $i \in [p]$, $A_{ij} \in \reals^{m_j \times n_i}$ for $(i,j) \in [p] \times [q]$ with full row rank, $\alpha \in \reals^q_+$, $\beta \in \reals^p_+$, and $\rho > 0$.
Its unique global optimum \cite{LNY22} computes the optimal constant for a family of inequalities in functional analysis and probability theory, including the Brascamp--Lieb inequality and the entropy power inequality \cite{AJN22}.
In particular, when the objective in \eqref{eq:3-prob} reduces to 
\[
f_\text{BL} (X) = -\log \det X + \sum_{j=1}^q \alpha_j \log \det (A_j X A_j^T)
\]
(with $\ones^T\alpha=1$) and the constraints disappear, the problem reduces to the computation of the classical Brascamp--Lieb constant \cite{SVY18}.
A recent paper \cite{WS22a} shows that the function $f_\text{BL}$ is geodesically convex on $\SV{n}{++}$ and DCA converges linearly to the global optimum when minimizing~$f_\text{BL}$.
However, in the general case, the objective in \eqref{eq:3-prob} may not be geodesically convex, and the analysis in \cite{WS22a} is not directly applicable.
But still, the general problem~\eqref{eq:3-prob} is a DC program in the form of~\eqref{eq:prob}, with
\[
g(X) = -\sum_{i=1}^p \beta_i \log \det X_i, \qquad h(X) = -\sum_{j=1}^q \alpha_j \log \det \Big(\sum_{i=1}^p A_{ij} X_i A_{ij}^T + \rho I_{m_j} \Big),
\] 
and $\cC = \{(X_1,\ldots,X_p) \mid 0 \leq X_i \leq C_i, \ i \in [p]\}$,
and we have the following result.
\begin{prop} \label{prop:3-conv}
Problem \eqref{eq:3-prob} satisfies \cref{assume}. 
DCA for \eqref{eq:3-prob} converges linearly to its minimizer.
\end{prop}

Apart from the convergence guarantee, the subproblem at $k$th DCA iteration is 
\begin{equation} \label{eq:3-subprob}
\begin{array}{ll}
\mini & -\sum\limits_{i=1}^p \beta_i \log \det U_i + \sum\limits_{i=1}^p \big\langle \sum\limits_{j=1}^q \alpha_j A_{ij}^T (\sum\limits_{i=1}^p A_{ij} X_i^{(k)} A_{ij}^T + \rho I_{m_j})^{-1} A_{ij}, U_i \big\rangle , \\
\st & 0 \preceq U_i \preceq C_i, \; \; i \in [p]
\end{array}
\end{equation} 
with variables $U_i \in \symm^{n_i}$, $i \in [p]$,
and can be solved efficiently by Bregman PDHG.
To see this, we first rewrite the subproblem \eqref{eq:3-subprob} in the form of \eqref{eq:prob-pdhg}: $F(\{U_i\}) = \sum_{i=1}^p F_i(U_i)$, $G(\{V_i\}) = \sum_{i=1}^p G_i(V_i)$,
\[
F_i(U_i) = -\beta_i \log \det U_i + \sum_{j=1}^q \alpha_j \Big\langle A_{ij}^T \Big(\sum_{i=1}^p A_{ij} X_i^{(k)} A_{ij}^T + \rho I_{m_j} \Big)^{-1} A_{ij}, U_i \Big\rangle + \delta_{\SV{n_i}{+}} (U_i),
\] 
$G_i = \delta_{\cC_i}$ is the indicator function of $\cC_i = \{X \in \symm^{n_i} \mid X \preceq C_i\}$, and the linear operator $\cA=\mathrm{Id}$ is the identity mapping.
Note that both functions $F$ and $G$ are separable in their respective arguments, while the original function $h$ is not separable in $X_1,\ldots,X_p$.
The $t$th iteration of Bregman PDHG takes the following updates for all $i \in [p]$
\begin{subequations} \label{eq:3-pdhg}
\begin{align} 
U_i^{(t+1)} &= \argmin_{U_i \in \SV{n_i}{+}} {\textstyle \Big( -\beta_i \log \det U_i + \sum\limits_{j=1}^q \alpha_j \big\langle A_{ij}^T (A_{ij} X_i^{(k)} A_{ij}^T + \rho I_{m_j})^{-1} A_{ij}, U_i \big\rangle} \nonumber \\
&\phantom{= \argmin_{U_i \in \SV{n_i}{+}}} \quad + \inprod{V_i^{(t)}}{U_i} + \tfrac{1}{\tau} d(U_i, U_i^{(t)}) \Big) \label{eq:3-pdhg-u} \\
V_i^{(t+1)} &= \prox_{\sigma \delta_{\cC_i}^\ast} (V_i^{(t)} + \sigma(2U_i^{(t+1)} - U_i^{(t)})), \label{eq:3-pdhg-v}
\end{align}
\end{subequations}
Again, each $V_i$-update is dominated by an eigen-decomposition.
For each $U_i$-update, we choose the kernel function $\phi_i(X) = -\log \det X$ with domain $\dom \phi_i = \SV{n_i}{++}$, and the resulting update~\eqref{eq:3-pdhg-u} has a closed-form expression.
Recall that the limit point $\{(\Uhat_i, \Vhat_i)\}$ of Bregman PDHG~\eqref{eq:3-pdhg} is the solution of the subproblem \eqref{eq:3-subprob}, and is also the $(k+1)$st DCA iteration, \ie, $X_i^{(k+1)} := \Uhat_i$, $i \in [p]$.
The $(U_i,V_i)$-update in Bregman PDHG \eqref{eq:3-pdhg} can be performed in parallel while the DCA iterates $\{X_i^{(k)}\}$ are coupled in the original function $h$ as well as its linearization in~\eqref{eq:3-subprob}.

\section{Numerical experiments} \label{sec:results}

In this section, we evaluate the performance of \algo~applied to the three DC programs in \cref{sec:probs}.
The numerical results verify that \algo~is able to solve the three classes of nonconvex programs with only a few number of DCA iterations and an efficient subproblem solver.

The experiments are carried out in MATLAB 2021b on a server with AMD Opteron Processor 6128 CPU and 32GB memory.
Since real data are not available for these problems, all the data in our experiments are synthetic, following \cite{WS22a,FFS23}.
Two classes of baselines are considered.
The nonlinear programming solver \fmincon~in Global Optimization Toolbox of Matlab is used to compare the optimality of the computed solutions.
For efficiency comparison, we also replace the subproblem solver with Euclidean PDHG, and \textsf{\small MOSEK} \cite{mosek}.
Regarding implementation details, the line search technique for Bregman PDHG \cite{JV22a,JV22b,CC22} is adopted for an adaptive choice of the stepsizes $\tau$, $\sigma$.

\begin{table}[t]
\begin{minipage}{0.31\textwidth}
  \centering
  \scalebox{0.80}{
    \begin{tabular}{@{}ccrr@{}} \toprule
      $n$ & algo & \multicolumn{1}{c}{obj. value} & runtime \\ \midrule
      \multirow{2}{*}{$100$} & \algo & $-4.2736$ & $6.82$  \\ & \fmincon & $-3.3742$ & $10.34$ \\ \midrule
      \multirow{2}{*}{$500$} & \algo & $-3.5771$ & $37.98$ \\ & \fmincon & $-1.2548$ & $126.76$ \\ \bottomrule
    \end{tabular}
  }
\end{minipage}%
\hspace{0.1in}
\begin{minipage}{0.31\textwidth}
  \centering
  \scalebox{0.80}{
    \begin{tabular}{@{}ccrr@{}} \toprule
      $n$ & algo & \multicolumn{1}{c}{obj. value} & runtime \\ \midrule
      \multirow{2}{*}{$100$} & \algo & $-4.2674$ & $8.63$  \\ & \fmincon & $-3.6485$ & $22.17$ \\ \midrule
      \multirow{2}{*}{$500$} & \algo & $-5.4356$ & $54.89$ \\ & \fmincon & $-1.3288$ & $262.13$ \\ \bottomrule
  \end{tabular}
  }
\end{minipage}%
\hspace{0.1in}
\begin{minipage}{0.31\textwidth}
  \centering
  \scalebox{0.80}{
    \begin{tabular}{@{}ccrr@{}} \toprule
      $n$ & algo & \multicolumn{1}{c}{obj. value} & runtime \\ \midrule
      \multirow{2}{*}{$100$} & \algo & $-4.4671$ & $13.28$ \\ & \fmincon & $-3.9126$ & $135.84$ \\ \midrule
      \multirow{2}{*}{$500$} & \algo & $-3.9523$ & $67.02$ \\ & \fmincon & $-2.5843$ & $561.57$ \\ \bottomrule
    \end{tabular}
  }
\end{minipage}%
\caption{Results of \algo~and \fmincon~for the problems \eqref{eq:1-prob} (left), \eqref{eq:2-prob} (middle), and \eqref{eq:3-prob} (right).
Objective function values and runtime (in sec.) are reported for a typical run of the algorithms.}
\label{tab:fmincon}
\end{table}

\begin{table}[t]
\centering
\scalebox{0.95}{
\begin{tabular}{@{}rlrrrr@{}} \toprule
  \multicolumn{1}{c}{$n$} & \multicolumn{1}{c}{algo} 
    & \multicolumn{1}{c}{\begin{tabular}[c]{@{}c@{}} avg. num. of \\ DCA iter. \end{tabular}} 
    & \multicolumn{1}{c}{\begin{tabular}[c]{@{}c@{}} avg. num. of \\ inner iter. \end{tabular}} 
    & \multicolumn{1}{c}{\begin{tabular}[c]{@{}c@{}} avg. runtime \\ (in sec.) \end{tabular}} 
    & \multicolumn{1}{c}{\begin{tabular}[c]{@{}c@{}} avg. runtime \\ per DCA iter. \end{tabular}} \\ \midrule
  \multirow{3}{*}{$500$} & \algo         & $9.5$ & $1735$ & $3.63 \times 10^2$ & $38.23$ \\
                          & \algo~(Euc.) & $9.5$ & $2046$ & $3.81 \times 10^2$ & $40.09$ \\
                          & \mosek       & $8.9$ & $76$   & $1.02 \times 10^3$ & $108.1$ \\ \midrule
  \multirow{3}{*}{$1000$} & \algo       &  $13.6$ & $1324$ & $1.73 \times 10^3$ & $127.2$ \\
                          & \algo~(Euc.) & $13.6$ & $1684$ & $2.20 \times 10^3$ & $162.4$ \\
                          & \mosek       & $13.2$ & $96$   & $9.87 \times 10^3$ & $726.3$ \\ \bottomrule
\end{tabular}
}
\caption{Results for problem \eqref{eq:1-prob}. Three algorithms are tested: \algo, \algo~with~Euc\-lidean PDHG, DCA with \textsf{\small MOSEK} as subproblem solver. The results are averaged over 10 synthetic datasets.}
\label{tab:result-1}
\end{table}

\cref{tab:fmincon} presents the results of \algo~and \fmincon.
Small-size problems are sufficient for comparison with \fmincon, and one typical result in the experiments is presented.
For all three classes of problems, \fmincon~might return a sub-optimal stationary point while \algo~is able to find a feasible point with a smaller objective value.
This is consistent with our analysis that for these problems, DCA converges to the unique global optimum (see \cref{thm:dca-conv}).

\cref{tab:result-1,tab:result-2,tab:result-3} present the numerical results for \eqref{eq:1-prob}, \eqref{eq:2-prob}, and \eqref{eq:3-prob}, respectively.
Three algorithms are tested: \algo, \algo~with Euclidean PDHG, and DCA with \textsf{\small MOSEK} as the subproblem solver.
The results are averaged over 10 experiments (with different synthetic data).
In all the experiments, DCA (the outer loop of \algo) converges within around 10 iterations, which supports the linear convergence result in \cref{thm:dca-conv}.
\textsf{\small DC-MOSEK} has a bit fewer DCA iterations because the interior-point solver \textsf{\small MOSEK} often returns a more accurate subproblem solution than first-order PDHG.
It is also observed that the customized subproblem solver in \algo~(Bregman PDHG) performs consistently better than the other two.
The average number of iterations in both PDHG methods is about one to two thousand, which is typical for these primal--dual first-order proximal methods.
(interior-point solver \textsf{\small MOSEK} takes much fewer iterations but has much more expensive per-iteration complexity.)
The per-iteration runtime of Bregman PDHG is smaller than that of its Euclidean variant, because the former does not need to solve a number of $n$ quadratic equations (as explained in \cref{sec:prob-1}).
Bregman PDHG takes fewer iterations (on average) than Euclidean PDHG, which might be explained by the fact that the Bregman distance can help build a more accurate local optimization model around the current iterate \cite{JV22b}.

\begin{table}[t]
\centering
\scalebox{0.95}{
\begin{tabular}{@{}rlrrrr@{}} \toprule
  \multicolumn{1}{c}{$n$} & \multicolumn{1}{c}{algo} 
    & \multicolumn{1}{c}{\begin{tabular}[c]{@{}c@{}} avg. num. of \\ DCA iter. \end{tabular}} 
    & \multicolumn{1}{c}{\begin{tabular}[c]{@{}c@{}} avg. num. of \\ inner iter. \end{tabular}} 
    & \multicolumn{1}{c}{\begin{tabular}[c]{@{}c@{}} avg. runtime \\ (in sec.) \end{tabular}} 
    & \multicolumn{1}{c}{\begin{tabular}[c]{@{}c@{}} avg. runtime \\ per DCA iter. \end{tabular}} \\ \midrule
  \multirow{3}{*}{$500$} & \algo        & $10.2$ & $1273$ & $5.71 \times 10^2$ & $56.07$ \\
                         & \algo~(Euc.) & $10.2$ & $1496$ & $5.63 \times 10^2$ & $75.83$ \\
                         & \mosek       & $9.8$  & $93$   & $2.32 \times 10^3$ & $225.1$ \\ \midrule
  \multirow{3}{*}{$1000$} & \algo       & $12.4$ & $1468$ & $3.50 \times 10^3$ & $281.9$ \\
                         & \algo~(Euc.) & $12.4$ & $1632$ & $4.08 \times 10^3$ & $313.3$ \\
                         & \mosek       & - & - & - & - \\ \bottomrule
\end{tabular}
}
\caption{Results for problem \eqref{eq:2-prob}. Three algorithms are tested: \algo, \algo~with~Euc\-lidean PDHG, DCA with \textsf{\small MOSEK} as the subproblem solver. The results are averaged over 10 synthetic datasets. `-' indicates the experiments do not finish in 6 hours.}
\label{tab:result-2}
\end{table}

\begin{table}
\centering
\scalebox{0.95}{
\begin{tabular}{@{}rlrrrr@{}} \toprule
  \multicolumn{1}{c}{$n$} & \multicolumn{1}{c}{algo} 
    & \multicolumn{1}{c}{\begin{tabular}[c]{@{}c@{}} avg. num. of \\ DCA iter. \end{tabular}} 
    & \multicolumn{1}{c}{\begin{tabular}[c]{@{}c@{}} avg. num. of \\ inner iter. \end{tabular}} 
    & \multicolumn{1}{c}{\begin{tabular}[c]{@{}c@{}} avg. runtime \\ (in sec.) \end{tabular}} 
    & \multicolumn{1}{c}{\begin{tabular}[c]{@{}c@{}} avg. runtime \\ per DCA iter. \end{tabular}} \\ \midrule
  \multirow{3}{*}{$500$} & \algo        & $14.7$ & $1157.9$ & $9.98 \times 10^2$ & $64.21$ \\
                         & \algo~(Euc.) & $14.7$ & $1297.5$ & $1.14 \times 10^3$ & $70.42$ \\
                         & \mosek       & $13.9$ & $85.2$   & $5.36 \times 10^4$ & $364.8$ \\ \midrule
  \multirow{3}{*}{$1000$} & \algo       & $14.2$ & $1048.7$ & $5.74 \times 10^3$ & $412.6$ \\
                         & \algo~(Euc.) & $14.2$ & $1362.6$ & $6.52 \times 10^3$ & $468.7$ \\
                         & \mosek & - & - & - & - \\ \bottomrule
\end{tabular}
}
\caption{Results for problem \eqref{eq:3-prob} ($p=q=3$, $n_i=n$ for $i \in [p]$). Three algorithms are tested: \algo, \algo~with Euclidean PDHG, DCA with \textsf{\small MOSEK} as the subproblem solver. The results are averaged over 10 synthetic datasets. `-' indicates the experiments do not finish in 6 hours.}
\label{tab:result-3}
\end{table}

\section{Conclusions} \label{sec:conclusion}

We present the \algo~framework for solving several classes of difference-of-convex (DC) programming problems that have a bounded, closed, and convex constraint set.
\algo~incorporates the difference-of-convex algorithm (DCA) in the outer iteration and a Bregman primal--dual proximal algorithm (Bregman PDHG) in the inner iteration.
In the theoretical aspect, the outer iteration DCA is shown to converge linearly to a global optimum under an extension of the classical Polyak--Łojasiewicz condition; from the computational perspective, Bregman PDHG solves the subproblem at each DCA iteration efficiently, with simple update rules and cheap per-iteration complexity.
As an application, we consider several classes of matrix optimization problems involving log-determinant functions.
For these problems, linear global convergence of DCA is established via the presented analysis, and Bregman PDHG efficiently solves the DCA subproblems by a suitable choice of Bregman distances.
Numerical experiments are conducted to verify the efficiency of \algo.
Future work are needed to identify infeasible and pathological DC programs in the aforementioned applications.

\newpage
\bibliography{dcprox}

\end{document}